\documentclass[11pt]{article}
\usepackage{anysize}
\marginsize{2.0cm}{2.0cm}{2.0 cm}{2.0 cm}
\usepackage{setspace}

\usepackage{latexsym,amsfonts,amsmath}
\usepackage{graphicx}
\usepackage{color}
\usepackage{epstopdf}
\usepackage[titletoc]{appendix}

\newtheorem{condition}{Condition}[section]{\bfseries}{\itshape}

\newtheorem{assumption}{Assumption}[section]{\bfseries}{\itshape}

\newtheorem{theorem}{Theorem}[section]{\bfseries}{\itshape}

{\bfseries}{\itshape}

{\bfseries}{\itshape}

{\bfseries}{\itshape}

\newtheorem{lemma}{Lemma}[section]{\bfseries}{\itshape}

{\bfseries}{\itshape}

\newtheorem{definition}{Definition}[section]{\bfseries}{\itshape}

{\bfseries}{\itshape}

\begin{document}

\title{Note on discounted continuous-time Markov decision processes with a lower bounding function}
\author{Xin Guo\thanks{Department of Mathematical Sciences, University of
Liverpool, Liverpool, L69 7ZL, U.K.. E-mail: X.Guo21@liv.ac.uk.},~ Alexey Piunovskiy\thanks{Department of Mathematical Sciences, University of
Liverpool, Liverpool, L69 7ZL, U.K.. E-mail: piunov@liv.ac.uk.}
~and Yi
Zhang \thanks{Corresponding author. Department of Mathematical Sciences, University of
Liverpool, Liverpool, L69 7ZL, U.K.. E-mail: yi.zhang@liv.ac.uk.}}
\date{}
\maketitle

\par\noindent{\bf Abstract:} In this paper, we consider the discounted continuous-time Markov decision process (CTMDP) with a lower bounding function. In this model, the negative part of each cost rate is bounded by the drift function, say $w$, whereas the positive part is allowed to be arbitrarily unbounded. Our focus is on the existence of a stationary optimal policy for the discounted CTMDP problems out of the more general class. Both constrained and unconstrained problems are considered. Our investigations are based on a useful transformation for nonhomogeneous Markov pure jump processes that has not yet been widely applied to the study of CTMDPs. This technique was not employed in previous literature, but it clarifies the roles of the imposed conditions in a rather transparent way. As a consequence, we withdraw and weaken several conditions commonly imposed in the literature.
\bigskip



\section{Introduction}
In this paper, we consider the discounted continuous-time Markov decision process (CTMDP) with a lower bounding function. In this model, the negative part of each cost rate is bounded by the drift function, say $w$, whereas the positive part is allowed to be arbitrarily unbounded. Our focus is on the existence of a stationary optimal policy for the discounted CTMDP problems out of the more general class. Both constrained and unconstrained problems are considered. Our investigations are based on a useful transformation for nonhomogeneous Markov pure jump processes that has not yet been widely applied to the study of CTMDPs.

Discounted CTMDPs have been studied intensively since the 1960s, with one of the first works being \cite{Rykov:1966}. Initially the theory is mainly developed for the finite state space models with bounded cost and transition rates. Later developments extend to models in a Borel state space with unbounded transition and cost rates, see e.g., \cite{Feinberg:2012OHL,Guo:2007,ABPZY:20102}. When the cost rates are unbounded from both above and below, a standard setup is to assume that there is a weight (or Lyapunov) function say $w$, bounding the growth of the absolute value of the cost rates and the transition rates in a suitable sense, so that the value function will be also bounded by this function $w$. Then the investigation is based on the applicability of Dynkin's formula to the class of $w$-bounded functions, for which some additional conditions must be also imposed. This line of reasoning was followed and demonstrated in the recent monographs \cite{Guo:2009,PrietoOHL:2012book} and the articles \cite{Blok:2015,ABPZY:20102}. If, as in the present paper, we only bound the growth of the negative part of each cost rate using the function $w$, which is thus called a lower bounding function, then the value function is in general not $w$-bounded. The approach based on the Dynkin's formula becomes less adequate.

On the other hand, thanks to the powerful Feinberg's reduction technique \cite{Feinberg:2004,Feinberg:2012OHL}, now it is well known that a discounted CTMDP problem is equivalent to a total undiscounted DTMDP (discrete-time Markov decision process) problem with the same action space. (By the way, Feinberg's reduction technique is different from and much more powerful than the uniformization technique, and its extension to the total undiscounted CTMDP problems is more delicate, see \cite{GuoZhang:2016,Piunovskiy:2015}.) This approach has been applied to studying the discounted CTMDP problem with arbitrarily unbounded transition rate and nonnegative cost rates, see \cite{Feinberg:2012OHL}. Nevertheless, the case, where the cost rates can take both positive and negative values, has never been treated with this approach, to the best of our knowledge. The reason is that when the transition rate is unbounded, the induced DTMDP is in general not absorbing, and the cost functions can take both positive and negative values. Without additional conditions, the studies for such DTMDPs, especially for constrained problems, are challenging and difficult, as demonstrated in \cite{FeinbergSonin:1996}, and are still underdeveloped, see e.g., \cite{Dufour:2013}.

Having said the above, discounted CTMDP problems with a lower bounding function have not been studied in the literature. The corresponding model in discounted discrete-time problems was treated in \cite{BauerleRieder:2011,Jaskiewicz:2011}, where the motivation for considering this type of cost functions was explained with applications to economics in \cite{Jaskiewicz:2011}. Note that they can be reduced to equivalent discounted problems with nonnegative cost functions using the method in \cite{vanderWal:1981}, see also \cite{Dufour:2016}. We shall demonstrate the continuous-time version of this technique. In \cite{BauerleRieder:2011}, this type of model was studied for a specific piecewise deterministic Markov decision process with jumps driven by a Poisson process, but following a different method based on the Young topology,  compared with the one here.

Our main contributions are as follows. Under conditions similar to those in \cite{Blok:2015}, we show the existence of a deterministic stationary (respectively, stationary) optimal policy for the unconstrained (respectively, constrained) discounted CTMDP problems with a lower bounding function. Our argument is based on a transformation for nonhomogeneous Markov pure jump processes, which, under some additional conditions, allows us to reduce the original problems to equivalent problems with nonnegative cost rates, so as for the Feinberg's reduction technique to apply. The roles of the additional conditions for this reduction are self-justified in a rather transparent way, as compared to the justification based on their relation to the Dynkin's formula, see \cite{Blok:2015}, which considers only the undiscounted problem with a $w$-bounded cost rate in a denumerable state space, and is restricted to stationary policies. With the better understanding of the roles of the conditions, even in the specific case, where the cost rates are bounded by the drift function $w$, we improve the existing results in \cite{Guo:2007,ABPZY:20102} by withdrawing and weakening several conditions assumed therein.

The rest of the paper is organized as follows. In Section \ref{GuoZhangSec2} we formulate the optimal control problems under consideration. The main statement is presented and proved in Section \ref{GuoPZhangSecMain}. The paper is finished with a conclusion in Section \ref{GuoPZhangConclusionSec}. Some auxiliary definitions and facts are included in the appendix.

\section{Model description and problem statement}\label{GuoZhangSec2}

The objective of this section is to describe briefly the
controlled process similarly to
\cite{Feinberg:2004,Feinberg:2012OHL,Kitaev:1995,ABPZY:20102}, and the associated
optimal control problem of interest in this paper.

In what follows, ${\cal{B}}(X)$ is
the Borel $\sigma$-algebra of the Borel space $X,$ $I$ stands for the indicator function, and $\delta_{\{x\}}(\cdot)$
is the Dirac measure concentrated on the singleton $\{x\}$. A measure is $\sigma$-additive and $[0,\infty]$-valued. Below, unless stated otherwise, the term of
measurability is always understood in the Borel sense. Throughout
this article, we adopt the conventions of
$
\frac{0}{0}:=0,~0\cdot\infty:=0,~\frac{1}{0}:=+\infty,~\infty-\infty:=\infty.
$

The primitives of a CTMDP are the following elements $\{S,
A, A(\cdot), q \},$ where $S$ is a nonempty Borel state space, $A$ is
a nonempty Borel action space, the ${\cal B}(A)$-valued multifunction $x\in S\rightarrow A(x)$ is, by assumption, with a measurable graph $\mathbb{K}:=\{(x,a)\in S\times A:~a\in A(x)\}$, and $q$
stands for a signed kernel $q(dy|x,a)$ on ${\cal{B}}(S)$ given
$(x,a)\in \mathbb{K}$ such that
$
\widetilde{q}(\Gamma|x,a):=q(\Gamma_S\setminus\{x\}|x,a)\ge 0
$
for all $\Gamma\in{\cal{B}}(S)$. Throughout this paper, we assume
that $q(\cdot|x,a)$ is conservative and stable, i.e.,
$
q(S|x,a)=0,~\bar{q}_x=\sup_{a\in A(x)}q_x(a)<\infty,
$
where $q_x(a):=-q(\{x\}|x,a).$ The signed kernel $q$ is often called the
transition rate. Below we assume that the set $\mathbb{K}$ contains the graph of some measurable mapping from $S$ to $A$.

Let us take the sample space $\Omega$ by adjoining to the
countable product space $S\times((0,\infty)\times S)^\infty$ the
sequences of the form
$(x_0,\theta_1,\dots,\theta_n,x_n,\infty,x_\infty,\infty,x_\infty,\dots),$
where $x_0,x_1,\dots,x_n$ belong to $S$,
$\theta_1,\dots,\theta_n$ belong to $(0,\infty),$ and
$x_{\infty}\notin S$ is the isolated point. We equip $\Omega$ with
its Borel $\sigma$-algebra $\cal F$.

Let $t_0(\omega):=0=:\theta_0,$ and for each $n\geq 0$, and each
element $\omega:=(x_0,\theta_1,x_1,\theta_2,\dots)\in \Omega$, let
$
t_n(\omega):=t_{n-1}(\omega)+\theta_n,
$
and
$
t_\infty(\omega):=\lim_{n\rightarrow\infty}t_n(\omega).
$
Obviously, $t_n(\omega)$ are measurable mappings on $(\Omega,{\cal
F})$. In what follows, we often omit the argument $\omega\in
\Omega$ from the presentation for simplicity. Also, we regard
$x_n$ and $\theta_{n+1}$ as the coordinate variables, and note
that the pairs $\{t_n,x_n\}$ form a marked point process with the
internal history $\{{\cal F}_t\}_{t\ge 0},$ i.e., the filtration
generated by $\{t_n,x_n\}$; see Chapter 4 of \cite{Kitaev:1995}
for greater details. The marked point process $\{t_n,x_n\}$
defines the stochastic process on $(\Omega,{\cal F})$ of interest
$\{\xi_t,t\ge 0\}$ by
\begin{eqnarray}\label{ZhangExponentialGZCTMDPdefxit}
\xi_t=\sum_{n\ge 0}I\{t_n\le t<t_{n+1}\}x_n+I\{t_\infty\le
t\}x_\infty.
\end{eqnarray}
Here we accept $0\cdot x:=0$ and $1\cdot x:=x$ for each $x\in S_\infty,$ and below we denote
$S_{\infty}:=S\bigcup\{x_\infty\}$.

\begin{definition}
\begin{itemize}
\item[(a)]
A policy $\pi$ for the CTMDP is a ${\cal P}(A)$-valued predictable process with respect to the internal history $\{{\cal F}_t\}$ so that, for each $\omega=(x_0,\theta_1,x_1,\theta_2,\dots)\in
\Omega$ and $t\in(0,\infty),$
\begin{eqnarray*}
\pi(da|\omega,t)&=&I\{t\ge t_\infty\}\delta_{a_\infty}(da)+
\sum_{n=0}^\infty I\{t_n< t\le
t_{n+1}\}\pi_{n}(da|x_0,\theta_1,\dots,\theta_n,x_n, t-t_n),
\end{eqnarray*}
where $a_\infty\notin A$ is some isolated point. Here, for each $n=0,1,2,\dots,$
$\pi_n(da|x_0,\theta_1,\dots,x_{n},s)$ is a stochastic kernel on
$A$ concentrated on $A(x_n)$ given $x_0\in S,~\theta_1\in (0,\infty),\dots,~x_{n}\in S,~s\in(0,\infty)$. We often identify a policy $\pi$ with the sequence of stochastic kernels $\{\pi_n\}_{n=0}^\infty.$

\item[(b)] A policy $\pi$ is called Markov if, for some stochastic kernel $\varphi$ on $A$ concentrated on $A(x)$ from $(x,t)\in S\times (0,\infty),$ one can write $\pi(da|\omega,t)=\varphi(da|\xi_{t-},t)$ whenever $t<t_\infty$  A Markov policy is identified with the underlying stochastic kernel $\varphi$.
\item [(c)] A policy
$\pi=\{\pi_n\}_{n=0}^\infty$  is called stationary if, with slight abuse of
notations, each of the stochastic kernels $\pi_n$ reads $
\pi_n(da|x_0,\theta_1,\dots,x_{n},s)=\pi(da|x_{n}). $ A
stationary policy is further called deterministic if $
\pi_n(da|x_0,\theta_1,\dots,x_{n},s)=\delta_{\{f(x_{n})\}}(da) $
for some measurable mapping $f$ from $S$ to $A$ such that $f(x)\in A(x)$ for each $x\in S$. We shall identify such a deterministic stationary policy by the underlying measurable mapping $f$.
\end{itemize}
\end{definition}
The class of all policies for the CTMDP is denoted by $\Pi,$ and the class of all Markov policies is $\Pi^M.$

Under a policy $\pi=\{\pi_n\}_{n=0}^\infty\in \Pi$, we define the
following predictable random measure $\nu^\pi$ on $S\times (0,\infty)$ by
\begin{eqnarray*}
\nu^\pi(dt, dy)&:=& \int_A \widetilde{q}(dy|\xi_{t-}(\omega),a)\pi(da|\omega,t)dt\nonumber\\
&=&\sum_{n\ge 0}\int_A
\widetilde{q}(dy|x_n,a)\pi_{n}(da|x_0,\theta_1,\dots, \theta_n,
x_n,t-t_n)I\{t_n< t\le t_{n+1}\}dt
\end{eqnarray*}
with $q_{x_\infty}(a_\infty)=q(dy|x_\infty,a_\infty):=0=:q_{x_\infty}(a)$ for each $a\in A.$ Then, given the initial distribution  $\gamma$, i.e., a probability measure on ${\cal B}(S)$,
there exists a unique probability measure ${P}^\pi_\gamma$ such
that
\begin{eqnarray*}
{P}_{\gamma}^\pi(x_0\in dx)=\gamma(dx),
\end{eqnarray*}
and with
respect to $P_\gamma^\pi,$ $\nu^\pi$ is the dual predictable
projection of the random measure associated with the marked point
process $\{t_n,x_n\}$; see \cite{Jacod:1975,Kitaev:1995}. Below,
when $\gamma$ is a Dirac measure concentrated at $x\in S,$
we use the denotation ${}{P}_x^\pi.$ Expectations with respect to
${}{P}_\gamma^\pi$ and ${}{P}_x^\pi$ are denoted as
${}{E}_{\gamma}^\pi$ and ${}{E}_{x}^\pi,$ respectively.

According to \cite{Jacod:1975}, the conditional distribution of $(\theta_{n+1},x_{n+1})$ with the condition on $x_0,\theta_1, \dots,\theta_{n},x_n$ is given on $\{\omega:x_n(\omega)\in S\}$ by
\begin{eqnarray*}
&&P_\gamma^\pi(\theta_{n+1}\in \Gamma_1,~x_{n+1}\in \Gamma_2|x_0,\theta_1,x_1,\dots,\theta_{n},x_n)\\
&=&\int_{\Gamma_1}e^{-\int_0^t \int_A q_{x_n}(a)\pi_n(da|x_0,\theta_1,\dots,\theta_n,x_n,s)ds}\int_{A}\widetilde{q}(\Gamma_2|x_n,a)\pi_n(da|x_0,\theta_1,\dots,\theta_n,x_n,t)dt,\\
&&~\forall~\Gamma_1\in{\cal B}((0,\infty)),~\Gamma_2\in{\cal B}(S);\\
&&P_\gamma^\pi(\theta_{n+1}=\infty,~x_{n+1}=x_\infty|x_0,\theta_1,x_1,\dots,\theta_{n},x_n)=e^{-\int_0^\infty  \int_A q_{x_n}(a)\pi_n(da|x_0,\theta_1,\dots,\theta_n,x_n,s)ds},
\end{eqnarray*}
and given on $\{\omega:x_n(\omega)=x_\infty\}$ by
\begin{eqnarray*}
P_\gamma^\pi(\theta_{n+1}=\infty,~x_{n+1}=x_\infty|x_0,\theta_1,x_1,\dots,\theta_{n},x_n)=1.
\end{eqnarray*}

Let $\infty>\alpha>0$ be a fixed discount factor. For each $j=0,1,\dots,N,$ with $N\ge 1$ being a fixed integer, let $c_j$ be a $(-\infty,\infty]$-valued measurable function on $\mathbb{K}$, representing a cost rate, and $d_j$ be a fixed finite constant, representing a corresponding constraint. We shall consider the following unconstrained and constrained $\alpha$-discounted optimal control problems for the CTMDP $\{S,A,A(\cdot),q\}$, respectively:
\begin{eqnarray}\label{ZyExponChap1}
\mbox{Minimize over $\pi\in \Pi$:}&&
E_x^\pi\left[\int_0^\infty  e^{-\alpha t}\int_A
c_0(\xi_t,a)\pi(da|\omega,t)dt\right],~x\in S,
\end{eqnarray}
and
\begin{eqnarray}\label{ZyNote01}
\mbox{Minimize over $\pi\in \Pi$:}&&
E_x^\pi\left[\int_0^\infty e^{-\alpha t}\int_A
c_0(\xi_t,a)\pi(da|\omega,t)dt\right]\nonumber\\
\mbox{such that}&& E_x^\pi\left[\int_0^\infty e^{-\alpha t}\int_A
c_j(\xi_t,a)\pi(da|\omega,t)dt\right]\le d_j,~j=1,2,\dots,N.
\end{eqnarray}
Here and below, we put
\begin{eqnarray}\label{ZYNote022}
c(x_\infty,a):=0,~\forall~a\in A\bigcup\{a_\infty\}.
\end{eqnarray}
The conditions we impose below will ensure that the performance measures in the above two problems are well defined, though not necessarily finite.

A policy $\pi^\ast$ is called optimal for the unconstrained problem (\ref{ZyExponChap1}) if
\begin{eqnarray*}
E_x^{\pi^\ast}\left[\int_0^\infty  e^{-\alpha t}\int_A
c_0(\xi_t,a)\pi^{\ast}(da|\omega,t)dt\right]=\inf_{\pi\in \Pi}E_x^\pi\left[\int_0^\infty  e^{-\alpha t}\int_A
c_0(\xi_t,a)\pi(da|\omega,t)dt\right],~\forall~x\in S.
\end{eqnarray*}

A policy $\pi$ is called feasible for the constrained problem (\ref{ZyNote01}) if
it satisfies all the inequalities therein. A feasible policy $\pi$ for problem (\ref{ZyNote01}) is said to be of a finite value if
 \begin{eqnarray*}
 -\infty<E_x^\pi\left[\int_0^\infty e^{-\alpha t}\int_A
c_0(\xi_t,a)\pi(da|\omega,t)dt\right]<\infty.
\end{eqnarray*}
 A policy $\pi^\ast$ is said to be optimal for problem (\ref{ZyNote01}) if it is feasible and satisfies
\begin{eqnarray*}
E_x^{\pi^\ast}\left[\int_0^\infty e^{-\alpha t}\int_A
c_0(\xi_t,a)\pi^\ast(da|\omega,t)dt\right]\le E_x^\pi\left[\int_0^\infty e^{-\alpha t}\int_A
c_0(\xi_t,a)\pi(da|\omega,t)dt\right]
\end{eqnarray*}
for each feasible policy $\pi$.

Note that the optimality of a feasible policy for the constrained problem (\ref{ZyNote01}) is for the fixed initial state $x\in S$. Here, we did not consider the more general case of a fixed initial distribution just for brevity and readability. The case of a fixed initial distribution $\gamma$ can be similarly treated with additional conditions regarding $\gamma$.

We would like to allow the possibility of cost rates unbounded from both above and below. We consider the following set of conditions to guarantee that the  performance measures in problems (\ref{ZyExponChap1}) and (\ref{ZyNote01}) are well defined.

\begin{condition}\label{ZYNoteCondition01}
There exists a $[1,\infty)$-valued measurable function $w$ on $S$ such that
\begin{itemize}
\item[(a)] for some finite constant $0<\rho<\alpha$,
\begin{eqnarray*}
\int_S w(y)q(dy|x,a)\le \rho w(x),~\forall~(x,a)\in\mathbb{K};
\end{eqnarray*}
\item[(b)] for some finite constant $L>0$,
\begin{eqnarray*}
c_i^-(x,a)\le L w(x),~\forall~(x,a)\in\mathbb{K},~i=0,1,\dots,N.
\end{eqnarray*}
Here, for each $i=0,1,\dots,N,$ $c_i^-$ is the negative part of the function $c_i$.
\end{itemize}
\end{condition}
Below, we accept that $w(x_\infty):=0$.
The cost rates satisfying part (b) of the above condition are said to be with the lower bounding function $w$; c.f. p.251 of \cite{BauerleRieder:2011} for a related definition for piecewise deterministic Markov decision processes.

\begin{lemma} Suppose Condition \ref{ZYNoteCondition01} is satisfied. Let a policy $\pi$ be arbitrarily fixed.
Then
\begin{eqnarray*}
E_x^\pi\left[\int_0^\infty  e^{-\alpha t}w(\xi_t)dt\right]<\infty,~\forall~x\in S.
\end{eqnarray*}
In particular, for each $x\in S,$ the integrals $E_x^\pi\left[\int_0^\infty  e^{-\alpha t}\int_A
c_i(\xi_t,a)\pi(da|\omega,t)dt\right],$ $i=0,1,\dots,N,$ are well defined.
\end{lemma}
\par\noindent\textit{Proof.} This follows from Lemma 2 of \cite{PiunovskiyZhang:2011Arxiv} and (\ref{ZYNote022}). $\hfill\Box$\bigskip

\begin{assumption}
Throughout this paper, unless stated otherwise, Condition \ref{ZYNoteCondition01} is assumed to hold automatically, without specific reference.
\end{assumption}

\section{Main statement and its proof}\label{GuoPZhangSecMain}
\subsection{Conditions, statements and comments}
\begin{condition}\label{ZhangNoteCondition02}
There exist a $(0,\infty)$-valued measurable function $w'$ on $S$ and a monotone nondecreasing sequence of measurable subsets $\{V_m\}_{m=1}^\infty\subseteq {\cal B}(S)$ such that the following hold.
\begin{itemize}
\item[(a)] $V_m\uparrow S$ as $m\rightarrow \infty.$
\item[(b)] $\sup_{x\in V_m}\overline{q}_x<\infty$ for each $m=1,2,\dots.$
\item[(c)] For some constant $\rho'\in(0,\infty),$
\begin{eqnarray*}
\int_S w'(y)q(dy|x,a)\le \rho' w'(x),~\forall~x\in S,~a\in A(x).
\end{eqnarray*}
\item[(d)] $\inf_{x\in S\setminus V_m}\frac{w'(x)}{w(x)}\rightarrow \infty$ as $m\rightarrow \infty,$ where the function $w$ comes from Condition \ref{ZYNoteCondition01}.
\end{itemize}
\end{condition}

Let a $[0,\infty)$-valued function $v$ on $S$ be fixed. A function $g$ on $S$ is called $v$-bounded if $||g||_v:=\sup_{x\in S}\frac{|g(x)|}{v(x)}<\infty;$ here the convention of $0/0=0$ is in use.
\begin{condition}\label{GZhangNoteCompactContinuity}
\begin{itemize}
\item[(a)] The multifunction $x\in S\rightarrow A(x)\in {\cal B}(A)$ is compact-valued and upper semicontinuous.
\item[(b)] For each $w$-bounded continuous function $g$ on $S$, $(x,a)\in \mathbb{K}\rightarrow \int_S g(y)\tilde{q}(dy|x,a)$ is continuous. Here and below the function $w$ is from Condition \ref{ZYNoteCondition01}.
\item[(b)] The function $w$ is continuous on $S$, and the functions $c_i$ are lower semicontinuous on $\mathbb{K}.$
\end{itemize}
\end{condition}

The next condition is for constrained problem only.
\begin{condition}\label{GuoZhangNoteConditionConstrained}
There exists a feasible policy for problem (\ref{ZyNote01}) with a finite value.
\end{condition}

The main statement of this paper is the following one.
\begin{theorem}\label{GuoXPZhangTheom1}
Suppose Conditions \ref{ZYNoteCondition01}, \ref{ZhangNoteCondition02} and \ref{GZhangNoteCompactContinuity} are satisfied. Then the following assertions hold.
\begin{itemize}
\item[(a)] There exists a deterministic stationary optimal policy for the unconstrained problem (\ref{ZyExponChap1}).
\item[(b)] If Condition \ref{GuoZhangNoteConditionConstrained} is also satisfied, then there exists a stationary optimal policy for the constrained problem (\ref{ZyNote01}).
\end{itemize}
\end{theorem}

In the previous literature, general discounted CTMDPs have not been considered when the cost rates were bounded below by a lower bounding function, and arbitrarily unbounded from the above, although for specific piecewise deterministic Markov decision processes with jumps driven by a Poisson process, this was considered in \cite{BauerleRieder:2011} following a different method. Discrete-time problems with a lower bounding function were considered in \cite{BauerleRieder:2011,Jaskiewicz:2011}, and in latter reference, the motivation for considering such cost functions was explained with their applications to economics. For discounted DTMDP problems, the treatment in \cite{BauerleRieder:2011,Jaskiewicz:2011} was direct. But it is possible to reduce this to equivalent problems with nonnegative cost functions, using the technique in p.101 of \cite{vanderWal:1981}, see also \cite{Dufour:2016} and p.79 of \cite{Altman:1999}. The proof of Theorem \ref{GuoXPZhangTheom1} will be based on a similar technique for CTMDPs, which, to the best of our knowledge, has not been widely applied to CTMDPs.

For the more restrictive case, where the cost rates are $w$-bounded, with $w$ coming from Condition \ref{ZYNoteCondition01}, Theorem \ref{GuoXPZhangTheom1}(a) was obtained in \cite{Blok:2015} under essentially equivalent conditions for discounted CTMDPs in a denumerable state space but restricted to the class of stationary policies. Our result here formally shows that it is without loss of generality to be restricted to this narrower class of policies under the imposed conditions. Otherwise, this sufficiency result seems not to follow from other known results in the relevant literature. The approach in \cite{Blok:2015} was directly based on the application of the Dynkin's forumla, and is different from ours. When the cost rates are only lower $w$-bounded, the value function is in general not $w$-bounded. Since under the conditions in \cite{Blok:2015} and here, Dynkin's formula is only applicable to the class of $w$-bounded functions, the treatment in \cite{Blok:2015} does not directly apply to the general case dealt with here.

Also when the cost rates are $w$-bounded, Theorem \ref{GuoXPZhangTheom1}(b) was obtained in e.g., \cite{ABPZY:20102} but under stronger conditions. We include them here for ease of reference.

Instead of Condition \ref{ZhangNoteCondition02}, the following condition was imposed in \cite{ABPZY:20102}.
\begin{condition}\label{GHOxINzHANGcONDITION}
There exists a $(0,\infty)$-valued measurable function $\widetilde{w}'$ on $S$ such that the following hold.
\begin{itemize}
\item[(a)] For some constant $\widetilde{L}'\in(0,\infty)$, $\overline{q}_x\le \widetilde{L}'\widetilde{w}'(x)$ for each $x\in S.$
\item[(b)] For some constant $\widetilde{\rho}'\in(0,\infty)$, $\int_S \widetilde{w}'(y)q(dy|x,a)\le \widetilde{\rho}' \widetilde{w}'(x)$ for each $(x,a)\in \mathbb{K}.$
\item[(c)] For some constant $\widetilde{L}\in(0,\infty)$, $(\overline{q}_x+1)w(x)\le \widetilde{L}\widetilde{w}'(x)$ for each $x\in S,$ where the function $w$ comes from Condition \ref{ZYNoteCondition01}.
\end{itemize}
\end{condition}

It is easy to see that, if the above condition is satisfied, then so is Condition \ref{ZhangNoteCondition02} with $w'=\widetilde{w}'+1$, $\rho'=\widetilde{\rho}'$, $V_m=\left\{x\in S:~\frac{\widetilde{w}'(x)+1}{w(x)}\le m\right\}$ for each $m=1,2,\dots.$

Furthermore, under Conditions \ref{ZYNoteCondition01}, \ref{ZhangNoteCondition02} and \ref{GuoZhangNoteConditionConstrained}, in addition to Condition \ref{GZhangNoteCompactContinuity}, it was also assumed in \cite{ABPZY:20102} that the function $\frac{\widetilde{w}'}{w}$ is a moment function on $\mathbb{K},$ see Definition E.7 of \cite{Hernandez-Lerma:1996}, in order to apply the Prokhorov theorem in their proof, see Proposition E.8 and Theorem E.6 of \cite{Hernandez-Lerma:1996}. This is not needed here. The investigations in \cite{ABPZY:20102} are largely based on the Dynkin's formula, and do not handle the more general cost rates considered here.

The rest of this section proves Theorem \ref{GuoXPZhangTheom1}. In the way, we comment and clarify the roles of the imposed conditions, and present the auxiliary statements.
\bigskip

\subsection{Proof of the main statement}
In this subsection, we present the proof of Theorem \ref{GuoXPZhangTheom1}, by combining several lemmas. To make the argument as transparent as possible, we proceed our proof in such a way that a lemma is presented only in the place, where it is needed in our proof, instead of collecting them altogether upfront.
\bigskip

\par\noindent\textit{\underline{Proof of Theorem \ref{GuoXPZhangTheom1}}.}
The following statement is a consequence of Theorem 4.2 of \cite{FeinbergShiyaevGood:2013}, and is the starting point of our reasoning.
\begin{lemma}
For each initial state $x\in S$ and policy $\pi,$ there exists a Markov policy $\varphi$ such that
\begin{eqnarray*}
E_x^\pi\left[\int_0^\infty  e^{-\alpha t}\int_A
f(\xi_t,a)\pi(da|\omega,t)dt\right]=E_x^\varphi\left[\int_0^\infty  e^{-\alpha t}\int_A
f(\xi_t,a)\varphi(da|\xi_t,t)dt\right]
\end{eqnarray*}
for each $[0,\infty]$-valued measurable function $f$ on $\mathbb{K}$.
\end{lemma}

The above lemma implies that without loss of generality, one can be restricted to the class of Markov policies for problems (\ref{ZyExponChap1}) and (\ref{ZyNote01}), i.e., if one obtains an optimal policy out of the class of Markov policies for problem (\ref{ZyExponChap1}) (or (\ref{ZyNote01})), then that policy is optimal for problem (\ref{ZyExponChap1}) (or (\ref{ZyNote01})) out of the general class.

We recall some definitions related to the process $\{\xi_t,t\ge 0\}$ under a Markov policy $\varphi$. Let us consider the signed kernel on $S$ from $S\times[0,\infty)$ defined by
\begin{eqnarray*}
q_\varphi(dy|x,t):=\int_A q(dy|x,a)\varphi(da|x,t),~\forall~x\in S,~t\in [0,\infty).
\end{eqnarray*}
Then $q_\varphi$ is a conservative and stable $Q$-function in the sense of \cite{FeinbergShiryayev:2014}, see p.262 therein. For the ease of reference, we recall some relevant definitions and facts about $Q$-functions in the appendix.

According to Theorem 2.2 of \cite{FeinbergShiryayev:2014}, under a Markov policy, say $\varphi$, the process $\{\xi_t,t\ge 0\}$ is a Markov pure jump process on $\{\Omega,{\cal F},\{{\cal F}_t \},P^\varphi\}$, that is,
  for each $s,t\in [0,\infty)$,
\begin{eqnarray*}
 P^\varphi(\xi_{t+s}\in \Gamma|{\cal F}_t)= P^\varphi(\xi_{t+s}\in \Gamma| \xi_t),~\forall~ \Gamma\in {\cal B}(X_\infty);
\end{eqnarray*}
and each trajectory of $\{\xi_t;t\ge 0\}$ is piecewise constant and right-continuous, such that for each $t\in[0,t_\infty)$, there are finitely many discontinuity points on the interval $[0,t]$. See Definition 1 in Chapter III of \cite{Gihman:1975}.
Here and below, we omit the subscript in $P^\varphi_\gamma$, whenever the initial distribution $\gamma$ is irrelevant. Furthermore, by Theorem 2.2 of \cite{FeinbergShiryayev:2014}, $p_{q_\varphi}$ defined by  (\ref{ZY2015:01}) with $q$ being replaced by $q_\varphi$ is the transition function corresponding to the process $\{\xi_t,t\ge 0\}$, i.e., for each $s\le t$, on $\{s<t_\infty\}$,
\begin{eqnarray*}
P^\varphi(\xi_t\in \Gamma|{\cal F}_s)=p_{q_\varphi}(s,\xi_s,t,\Gamma), ~\forall~\Gamma\in {\cal B}(S).
\end{eqnarray*}
(C.f. p.1397 of \cite{Kuznetsov:1984}.)
Consequently, for each Markov policy $\varphi$,
 \begin{eqnarray*}
 E_x^\varphi\left[\int_0^\infty  e^{-\alpha t}\int_A
c_i(\xi_t,a)\varphi(da|\xi_t,t)dt\right]=\int_0^\infty \int_{S} e^{-\alpha t}\int_A c_i(y,a)\varphi(da|y,t)p_{q_\varphi}(0,x,t,dy)dt,~\forall~x\in S
\end{eqnarray*}
for each  $i=0,1,\dots,N$.

Given the $Q$-function $q_\varphi$ on $S$ induced by a Markov policy $\varphi,$ let us introduce the $w$-transformed $Q$-function $q_\varphi^w$ on $S_\delta$ defined as follows.

Let \begin{eqnarray*}
S_\delta:=S\bigcup \{\delta\}
\end{eqnarray*}
with $\delta\notin S$ being an isolated point concerning the topology of $S_\delta$ that satisfies $\delta\ne x_\infty$. The $w$-transformed (stable conservative) $Q$-function $q_\varphi^w$ on $S_\delta$ is defined by
\begin{eqnarray}\label{ZY2015:08}
&&q_\varphi^w(\Gamma|x,s):=\left\{
\begin{array}{ll}
      \frac{\int_{\Gamma}w(y)q_\varphi(dy|x,s)}{w(x)}, & \mbox{if~}  x\in S,~ \Gamma\in {\cal B}(S),~x\notin \Gamma; \\
     \rho-\frac{\int_{S}w(y)q_\varphi(dy|x,s)}{w(x)}, & \mbox{if~} x\in S,~ \Gamma=\{\delta\};\\
     0, & \mbox{if~} x=\delta,~ \Gamma=S_\delta.
\end{array}
\right.
\end{eqnarray}
for each $s\in[0,\infty);$ and
\begin{eqnarray*}
{q_{\varphi}^w}_x(s):=\rho+{q_\varphi}_x(s),~\forall~s\in [0,\infty).
\end{eqnarray*}
Here, ${q_\varphi}_x(s)=-q_{\varphi}(S\setminus\{x\}|x,s)$; see the appendix for more definitions and relevant notations concerning a $Q$-function. For (uncontrolled) homogeneous continuous-time Markov chains, this transformation was considered in e.g., \cite{Anderson:1991,Spieksma:2015,Spieksma:2012}. But it has not been widely applied to the study of CTMDPs.

\begin{lemma}\label{ZY2015:09}
Let a Markov policy $\varphi$ be fixed. For each $x\in S,$ $s,t\in[0,\infty)$, $s\le t$ and $\Gamma\in {\cal B}(S)$, the following relation holds;
\begin{eqnarray*}
p_{q_\varphi^w}(s,x,t,\Gamma)=\frac{e^{-\rho(t-s)}}{w(x)}\int_\Gamma w(y) p_{q_\varphi}(s,x,t,dy).
\end{eqnarray*}
\end{lemma}
 \par\noindent\textit{Proof.} See Lemma A.3 of \cite{Zhang:2016}. $\hfill\Box$\bigskip

By Lemma \ref{ZY2015:09}, we see that for each $i=0,1,\dots,N,$
\begin{eqnarray*}
&&w(x)\int_0^\infty \int_S p_{q_\varphi^w}(0,x,t,dy)\int_A \frac{c_i(y,a)}{w(y)}\varphi(da|y,t) e^{-(\alpha-\rho)t}dt\\
&=&\int_0^\infty \int_S \int_A c_i(y,a)\varphi(da|y,t)e^{-\alpha t}p_{q_\varphi}(0,x,t,dy)dt,~\forall~x\in S.
\end{eqnarray*}
Hence,  problem (\ref{ZyExponChap1}) is equivalent to
\begin{eqnarray}\label{ZyNoteNOte011}
\mbox{Minimize over $\varphi\in \Pi^M$:}&&
\int_0^\infty \int_S p_{q_\varphi^w}(0,x,t,dy)\int_A \frac{c_0(y,a)}{w(y)}\varphi(da|y,t) e^{-(\alpha-\rho)t}dt,~x\in S,
\end{eqnarray}
and problem (\ref{ZyNote01}) is equivalent to
\begin{eqnarray}\label{ZyNoteNOte012}
\mbox{Minimize over $\varphi\in \Pi^M$:}&&
\int_0^\infty \int_S p_{q_\varphi^w}(0,x,t,dy)\int_A \frac{c_i(y,a)}{w(y)}\varphi(da|y,t) e^{-(\alpha-\rho)t}dt\nonumber\\
\mbox{such that}&& \int_0^\infty \int_S p_{q_\varphi^w}(0,x,t,dy)\int_A \frac{c_j(y,a)}{w(y)}\varphi(da|y,t) e^{-(\alpha-\rho)t}dt\le \frac{d_j}{w(x)},\nonumber\\
&&~j=1,2,\dots,N.
\end{eqnarray}
Thus, one can consider the $w$-transformed CTMDP $\{S_\delta,A\bigcup\{a_\infty\} ,A_\delta(\cdot),q^w\}$, where
$A_\delta(\delta):=\{a_\infty\}$, and $A_\delta(x):=A(x)$ for each $x\in S,$ while the transition rate $q^w$ is defined by
\begin{eqnarray*}
&&q^w(\Gamma|x,a)=\left\{
\begin{array}{ll}
      \frac{\int_{\Gamma}w(y)q(dy|x,a)}{w(x)}, & \mbox{if~}  x\in S,~ \Gamma\in {\cal B}(S),~x\notin \Gamma; \\
     \rho-\frac{\int_{S}w(y)q(dy|x,a)}{w(x)}, & \mbox{if~} x\in S,~ \Gamma=\{\delta\};\\
     0, & \mbox{if~} x=\delta,~ \Gamma=S_\delta.
\end{array}
\right.
\end{eqnarray*}
for each $x\in S_\delta$ and $a\in A_\delta(x)$; and
\begin{eqnarray*}
q^w_x(a):=\rho+q_x(a),~\forall~x\in S,~a\in A_\delta(x).
\end{eqnarray*}
The requirement of $\alpha>\rho$ in Condition \ref{ZYNoteCondition01}(a) is needed so that problems (\ref{ZyNoteNOte011}) and (\ref{ZyNoteNOte012}) are legitimate $(\alpha-\rho)$-discounted problems of the $w$-transformed CTMDP with the cost rates $c^w_i$ defined by
\begin{eqnarray*}
c^w_i(x,a):=\frac{c_i(x,a)}{w(x)}
\end{eqnarray*}
for each $x\in S,$ $a\in A(x)$; and
\begin{eqnarray*}
c^w_i(\delta,a_\infty):=0.
\end{eqnarray*}
According to the Feinberg's reduction technique for discounted CTMDPs, see \cite{Feinberg:2012OHL}, the CTMDP problems (\ref{ZyNoteNOte011}) and (\ref{ZyNoteNOte012}) can be reduced to equivalent total undiscounted problems for the DTMDP $\{S_\delta\bigcup\{x_\infty\}, A\bigcup\{a_\infty\}, A_\delta(\cdot),T\}$ with the cost functions $C_i$, where
the transition probability $T$ is defined by
\begin{eqnarray*}
T(\Gamma|x,a):=\frac{\int_\Gamma w(y)q(dy|x,a)}{(\alpha+q_x(a))w(x)}
\end{eqnarray*}
for each $\Gamma\in {\cal B}(S)$, $x\notin \Gamma,$ and $a\in A_\delta(x)$;
\begin{eqnarray*}
T(\{\delta\}|x,a):=\frac{\rho w(x)-\int_S w(y)q(dy|x,a)}{(\alpha+q_x(a))w(x)}
\end{eqnarray*}
for each $x\in S$ and $a\in A_\delta(x)$;
\begin{eqnarray*}
T(\{x_\infty\}|x,a):=\frac{\alpha-\rho}{\alpha+q_x(a)}
\end{eqnarray*}
for each $x\in S$ and $a\in A_\delta(x)$;
and  $T(\{x_\infty\}|x_\infty,a_\infty):=1=:T(\{x_\infty\}|\delta,a_\infty)$,
and the cost functions $C_i$ are defined by
\begin{eqnarray*}
C_i(x,a):=\frac{c_i(x,a)}{(\alpha+q_x(a))w(x)}
\end{eqnarray*}
for each $x\in S$ and $a\in A_\delta(x)$; and
\begin{eqnarray*}
C_i(\delta,a_\infty):=0=:C_i(x_\infty,a_\infty).
\end{eqnarray*}
More precisely, given the initial state $x\in S$, for each Markov policy $\varphi$ for the $w$-transformed CTMDP, there is a strategy $\sigma$ for the DTMDP $\{S_\delta\bigcup\{x_\infty\}, A\bigcup\{a_\infty\}, A_\delta(\cdot),T\}$ such that
\begin{eqnarray*}
\int_0^\infty \int_S p_{q_\varphi^w}(0,x,t,dy)\frac{c_i(y,a)}{w(y)} e^{-(\alpha-\rho)t}dt=\mathbb{E}_x^\sigma\left[\sum_{n=0}^\infty C_i(X_n,A_{n}) \right]
\end{eqnarray*}
for each $i=0,1,\dots,N$, and vice versa. Moreover, in the previous equality, if $\varphi$ is a deterministic stationary (respectively, stationary) policy, then $\sigma$ can be taken as a deterministic stationary (respectively, stationary) strategy for the DTMDP, and vice versa.  Here we use $\mathbb{E}_x^\sigma$ to denote the expectation taken with respect to the strategic measure of the DTMDP under the strategy $\sigma$, and $\{X_n\}$ and $\{A_{n}\}$ are the controlled and controlling processes in the DTMDP. The term ``strategy'' is reserved for the DTMDP to avoid the potential confusion with the corresponding notion for the CTMDP. We refer the reader to e.g., \cite{Hernandez-Lerma:1996,Piunovskiy:1997} for the standard description of a DTMDP.

Note that in general, the DTMDP $\{S_\delta\bigcup\{x_\infty\}, A\bigcup\{a_\infty\}, A_\delta(\cdot),T\}$ is not absorbing, and the cost function $C_i$ can take both positive and negative values. (This is the case e.g., if the original CTMDP is an uncontrolled pure birth process with $S=\{1,2,\dots\}$, and birth rate at the state $x\in S$ being $2x$, $\alpha=2,$ $\rho=1$ and $w(x)=1$ for each $x\in S$.) Compared to the absorbing model treated in \cite{Altman:1999,Feinberg:2012}, the theory for such a DTMDP model is technical and demanding, and, without additional assumptions, there is far less result concerning the existence of stationary strategies, which one can directly refer to, especially for the constrained problems, see \cite{Dufour:2013,FeinbergSonin:1996}.

On the other hand,  the functions $c^w_i$, $i=0,1,\dots,N,$ are bounded from below under Condition \ref{ZYNoteCondition01}(b). Let some common lower bound be $\underline{c}\le 0$. Let
\begin{eqnarray}\label{ZZZZZZ}
\widetilde{c}^w_i:=c_i^w-\underline{c}
\end{eqnarray}
for each $i=0,1,\dots,N.$ Then the functions $\widetilde{c}^w_i$ are all nonnegative. In order for problems  (\ref{ZyNoteNOte011}) and (\ref{ZyNoteNOte012}) to be equivalent to
\begin{eqnarray}\label{ZYEasyProblem1}
\mbox{Minimize over $\varphi\in \Pi^M$:}&&
\int_0^\infty \int_{S_\delta} p_{q_\varphi^w}(0,x,t,dy)\int_{A_\delta} \widetilde{c}^w_0(y,a)\varphi(da|y,t) e^{-(\alpha-\rho)t}dt,~x\in S,
\end{eqnarray}
and
\begin{eqnarray}\label{ZYEasyProblem2}
\mbox{Minimize over $\varphi\in \Pi^M$:}&&
\int_0^\infty \int_{S_\delta} p_{q_\varphi^w}(0,x,t,dy)\int_{A_\delta} \widetilde{c}_0^w(y)\varphi(da|y,t) e^{-(\alpha-\rho)t}dt\nonumber\\
\mbox{such that}&& \int_0^\infty \int_{S_\delta} p_{q_\varphi^w}(0,x,t,dy)\int_{A_\delta} \widetilde{c}_j^w(y)\varphi(da|y,t) e^{-(\alpha-\rho)t}dt\le \frac{d_j}{w(x)}-\frac{\underline{c}}{\alpha-\rho},\nonumber\\
&&~j=1,2,\dots,N,
\end{eqnarray}
respectively, we need the following relation to hold for each $\varphi\in \Pi^M$:
\begin{eqnarray}\label{ZhangNoteEquationooo}
p_{q_\varphi^w}(0,x,t,S_\delta)=1,~\forall~x\in S,~t\in[0,\infty).
\end{eqnarray}
Condition \ref{ZhangNoteCondition02} is precisely imposed for this purpose, as seen in the next statement. (An alternative justification of the role of Condition \ref{ZhangNoteCondition02} is that it validates the Dynkin's formula for the original CTMDP to a certain class of functions, see \cite{Blok:2015} for the homogeneous denumerable case. But the justification here is more transparent in our opinion.)

\begin{lemma} Let some Markov policy $\varphi$ be fixed. Suppose Condition \ref{ZYNoteCondition01}(a) and Condition \ref{ZhangNoteCondition02} are satisfied. Then (\ref{ZhangNoteEquationooo}) holds.
\end{lemma}
\par\noindent\textit{Proof.} According to Theorem \ref{ZY2015Theorem03}, for the statement it suffices to verify that Condition \ref{ZY2015Condition2} is satisfied.

Since the Markov policy $\varphi$ is fixed throughout this proof, we write $q_\varphi$ as $q$ for brevity.
Note that
\begin{eqnarray}\label{C3nnnnnnew1}
&&\int_{S}\frac{w'(y)}{w(y)}q^{w}(dy|x,s)=\int_{S}\frac{w'(y)}{w(y)}\frac{w(y)}{ w(x)}\widetilde{q}(dy|x,s)-(\rho+q_x(s))\frac{ w'(x)}{ w(x)} \nonumber\\
&=&\int_{S}\frac{ w'(y)}{ w(x)}\widetilde{q}(dy|x,s)-(\rho+q_x(s))\frac{ w'(x)}{ w(x)}\le (\rho'-\rho)\frac{w'(x)}{w(x)},~\forall~x\in S,~s\ge 0.
\end{eqnarray}
Consider the $[0,\infty)$-valued measurable function $\widetilde{w}$ on $[0,\infty)\times S_\delta$ defined for each $v\in[0,\infty)$ by  $\widetilde{w}(v,x)=\frac{w'(x)}{w(x)}$ if $x\in S$ and $\widetilde{w}(v,\delta)=0$. Then Condition \ref{ZY2015Condition2}, with $S$ and $q$ being replaced by $S_\delta$ and $q^{w}$,  is satisfied by the monotone nondecreasing sequence of measurable subsets $\{\widetilde{V}_n\}_{n=1}^\infty$ of $\mathbb{R}_+^0\times S_\delta$ defined by
$
\widetilde{V}_n=[0,\infty)\times V_n\bigcup\{\delta\}$ for each $n=1,2,\dots,
$
and the function $\widetilde{w}$ on $[0,\infty)\times S_\delta$ defined in the above. In greater detail, part (d) of the corresponding version of Condition \ref{ZY2015Condition2} is satisfied because, by (\ref{C3nnnnnnew1}),
\begin{eqnarray*}
&&\int_0^\infty \int_{S_\delta} \widetilde{w}(t+v,y) e^{-\rho' t -\int_{(0,t]} q^{w}_x(s+v)ds} \widetilde{{q^{w}}}(dy|x,t+v)dt\\
&\le& \int_0^\infty e^{-\rho' t -\int_0^t q^{w}_x(s+v)ds} \left(q_x(s)+\rho'\right)\widetilde{w}(v,x)
=\widetilde{w}(v,x),~\forall~x\in S,
\end{eqnarray*}
and the last inequality holds trivially when $x=\delta.$

Thus, by Theorem \ref{ZY2015Theorem03},  we see that relation (\ref{ZhangNoteEquationooo}) is satisfied, and the statement follows.  $\hfill\Box$
\bigskip

By the way, under Condition \ref{ZYNoteCondition01}(a), in certain models, Condition \ref{ZhangNoteCondition02} is also necessary for (\ref{ZhangNoteEquationooo}) to hold under certain policies; see \cite{Zhang:2016}. In the homogeneous denumerable case, this was first observed in \cite{Spieksma:2015}. For more concrete examples such as single birth processes, this necessity part was known earlier, see \cite{Chen:2015}.

As a result of the above lemma and the discussions above it, we see that under Condition \ref{ZYNoteCondition01} and Condition \ref{ZhangNoteCondition02}, one can reduce the $\alpha$-discounted problems (\ref{ZyExponChap1}) and (\ref{ZyNote01}) for the original CTMDP $\{S,A,A(\cdot),q\}$ to the $(\alpha-\rho)$-discounted problems (\ref{ZYEasyProblem1}) and (\ref{ZYEasyProblem2}) for the CTMDP $\{S_\delta, A_\delta, A_\delta(\cdot),q^w\}$ with nonnegative cost rates. Furthermore, according to the Feinberg's reduction technique \cite{Feinberg:2012OHL}, which was also sketched in the above, problems (\ref{ZYEasyProblem1}) and (\ref{ZYEasyProblem2})  can be reduced to
\begin{eqnarray}\label{GuoZyNoteReducedProbl1}
\mbox{Minimize over $\sigma$}&&
\mathbb{E}_x^\sigma\left[\sum_{n=0}^\infty \widetilde{C}_0(X_n,A_{n}) \right],~x\in S,
\end{eqnarray}
and
\begin{eqnarray}\label{GuoZyNoteReducedProbl2}
\mbox{Minimize over $\sigma$:}&&
\mathbb{E}_x^\sigma\left[\sum_{n=0}^\infty \widetilde{C}_0(X_n,A_{n}) \right]\nonumber\\
\mbox{such that}&& \mathbb{E}_x^\sigma\left[\sum_{n=0}^\infty \widetilde{C}_j(X_n,A_{n}) \right]\le \frac{d_j}{w(x)}-\frac{\underline{c}}{\alpha-\rho},\nonumber\\
&&~j=1,2,\dots,N,
\end{eqnarray}
respectively, for the DTMDP $\{S_\delta\bigcup\{x_\infty\}, A\bigcup\{a_\infty\}, A_\delta(\cdot),T\}$ defined earlier. Here the cost functions $\widetilde{C}_i$ for the DTMDP are defined by
\begin{eqnarray*}
\widetilde{C}_i(x,a):=\frac{\widetilde{c}^w_i(x,a)}{(\alpha+q_x(a)) }\ge 0
\end{eqnarray*}
for each $x\in S_\delta$ and $a\in A_\delta(x)$; and
\begin{eqnarray*}
\widetilde{C}_i(x_\infty,a_\infty):=0,
\end{eqnarray*}
with the functions $\widetilde{c}_i^w$ being defined by (\ref{ZZZZZZ}). Note that the cost functions $\widetilde{C}_i$ could be arbitrarily unbounded from above.

Finally, if Condition \ref{ZYNoteCondition01}, Condition \ref{ZhangNoteCondition02}, and Condition \ref{GZhangNoteCompactContinuity} are all satisfied, then it is easy to check that the DTMDP $\{S_\delta\bigcup\{x_\infty\}, A\bigcup\{a_\infty\}, A_\delta(\cdot),T\}$ with the nonnegative cost functions $\widetilde{C}_i$ is a semicontinuous  model, see \cite{BauerleRieder:2011,Dynkin:1979}, and it is a standard result that there exists an optimal deterministic stationary strategy for problem (\ref{GuoZyNoteReducedProbl1}). For the constrained problem (\ref{GuoZyNoteReducedProbl2}), under the extra Condition \ref{GuoZhangNoteConditionConstrained}, one can refer to Theorem 4.1 of \cite{Dufour:2012}, see also Theorem A.2 of \cite{CostDufour:2014}, for the existence of a stationary optimal strategy for (\ref{GuoZyNoteReducedProbl2}). Since these two DTMDP problems are equivalent to the original CTMDP problems, according to the Feinberg's reduction technique for discounted CTMDP problems as mentioned earlier, we immediately conclude the existence of an optimal deterministic stationary policy for the unconstrained CTMDP problem (\ref{ZyExponChap1}) and an optimal stationary policy for the constrained CTMDP problem (\ref{ZyNote01}). $\hfill\Box$
\bigskip

We finish this section with the following remark. In general, problems  (\ref{ZyNoteNOte011}) and (\ref{ZyNoteNOte012}) are not equivalent to
(\ref{ZYEasyProblem1}) and (\ref{ZYEasyProblem2}), respectively. According to \cite{Feinberg:2012OHL}, (\ref{ZYEasyProblem1}) is equivalent to the DTMDP problem $\{S_\delta\bigcup\{x_\infty\}, A\bigcup\{a_\infty\}, A_\delta(\cdot),T\}$ with the cost function $\widetilde{C}_0$. Suppose $\varphi^\ast$ is an optimal deterministic strategy for this DTMDP problem. Under Conditions  \ref{ZYNoteCondition01} and Condition \ref{GZhangNoteCompactContinuity}, if $V^\ast$ denotes the value function of this DTMDP problem, then such an optimal deterministic stationary strategy exists and can be obtained by taking the measurable selector providing the minimum in the following:
\begin{eqnarray*}
V^\ast(x)=\inf_{a\in A_\delta(x)}\left\{\widetilde{C}_0(x,a)+\int_{S_\delta}T(dy|x,a)V^\ast(y)\right\},~\forall~x\in S_\delta.
\end{eqnarray*}

We claim that $\varphi^\ast$ is also an optimal deterministic policy for the CTMDP problem (\ref{ZyNoteNOte011}), provided that (\ref{ZhangNoteEquationooo}) holds for this particular strategy $\varphi^\ast,$ i.e.,
\begin{eqnarray}\label{GuoXinPiunovPiunovZhang}
p_{q_{\varphi^\ast}^w}(0,x,t,S_\delta)=1,~\forall~x\in S,~t\in[0,\infty).
\end{eqnarray}
Indeed, since $\varphi^\ast$ is optimal for the DTMDP $\{S_\delta\bigcup\{x_\infty\}, A\bigcup\{a_\infty\}, A_\delta(\cdot),T\}$ with the cost function $\widetilde{C}_0,$ which is equivalent to  problem (\ref{ZYEasyProblem1}),
\begin{eqnarray*}
&&\inf_{\varphi\in \Pi^M}\left\{\int_0^\infty \int_{S_\delta} p_{q_\varphi^w}(0,x,t,dy)\int_{A_\delta} \widetilde{c}^w_0(y,a)\varphi(da|y,t)e^{-(\alpha-\rho)t}dt\right\}\\
&=&\int_0^\infty \int_{S_\delta} p_{q_{\varphi^\ast}^w}(0,x,t,dy)  \widetilde{c}^w_0(y,\varphi^\ast(y)) e^{-(\alpha-\rho)t}dt\\
&=&\int_0^\infty \int_{S} p_{q_{\varphi^\ast}^w}(0,x,t,dy) \frac{c_0(y,\varphi^\ast(y))}{w(y)}  e^{-(\alpha-\rho)t}dt-\frac{\underline{c}}{\alpha-\rho},~\forall~x\in S.
\end{eqnarray*}
Consider an arbitrarily fixed $\varphi\in \Pi^M.$ Then for each $x\in S,$
\begin{eqnarray*}
&&\int_0^\infty \int_{S} p_{q_{\varphi^\ast}^w}(0,x,t,dy) \frac{c_0(y,\varphi^\ast(y))}{w(y)}  e^{-(\alpha-\rho)t}dt-\frac{\underline{c}}{\alpha-\rho}\\
&\le&\int_0^\infty \int_{S_\delta} p_{q_\varphi^w}(0,x,t,dy)\int_{A_\delta} \widetilde{c}^w_0(y,a)\varphi(da|y,t)e^{-(\alpha-\rho)t}dt\\
&=&\int_0^\infty \int_{S} p_{q_\varphi^w}(0,x,t,dy)\int_{A} \frac{c_0(y,a)}{w(y)} \varphi(da|y,t)e^{-(\alpha-\rho)t}dt-\underline{c}\int_0^\infty p_{q_\varphi^w}(0,x,t,S_\delta)e^{-(\alpha-\rho)t}dt.
\end{eqnarray*}
Since $\underline{c}\le 0$, and $p_{q_\varphi^w}(0,x,t,S_\delta)\le 1$, it follows that
\begin{eqnarray*}
&&\int_0^\infty \int_{S} p_{q_{\varphi^\ast}^w}(0,x,t,dy) \frac{c_0(y,\varphi^\ast(y))}{w(y)}  e^{-(\alpha-\rho)t}dt\\
&\le&\int_0^\infty \int_{S} p_{q_\varphi^w}(0,x,t,dy)\int_{A} \frac{c_0(y,a)}{w(y)} \varphi(da|y,t)e^{-(\alpha-\rho)t}dt,~\forall~x\in S.
\end{eqnarray*}
Condition (\ref{GuoXinPiunovPiunovZhang}) can be checked using Theorem \ref{ZY2015Theorem03} in the appendix.  The similar reasoning also holds for the constrained problem. To avoid repetition, we omit the details.

\section{Conclusion}\label{GuoPZhangConclusionSec}
To sum up, we showed the existence of a deterministic stationary (respectively, stationary) optimal policy for the unconstrained (respectively, constrained) discounted CTMDP problems under rather weak conditions. The main feature in the model is that only the negative part of each cost rate is bounded by a drift function. Another contribution is that our arguments were based on a transformation for Markov pure jump processes, and this technique  had not been widely applied to the study of CTMDPs. On the other hand, exactly this technique allowed us to clarify the roles of all the imposed conditions in a transparent way. In this way, even in the specific case, where both the negative and positive parts of the cost rates are bounded by the drift function, we improved the existing results in the literature by withdrawing several and various conditions assumed therein.

\appendix
\section{Appendix}

A (Borel-measurable) signed kernel $q(dy|x,s)$ on ${\cal B}(S)$ from $S\times [0,\infty)$ is called a (conservative stable) $Q$-function on the Borel space $S$ if the following conditions are satisfied.
\begin{itemize}
\item[(a)] For each $s\ge 0$, $x\in S$ and $\Gamma\in {\cal B}(S)$ with $x\notin \Gamma,$
$
\infty>q(\Gamma|x,s)\ge 0.
$
\item[(b)]  For each $(x,s)\in S\times [0,\infty),$
$
q(S|x,s)=0.
$
\item[(c)] For each $x\in S,$
$
\sup_{s\in[0,\infty)}\left\{q(S\setminus \{x\}|x,s)\right\}<\infty.
$
\end{itemize}
For each $Q$-function $q$ on $S$, we put $\widetilde{q}(\Gamma|x,s):=q(\Gamma\setminus \{x\}|x,s)$, and $q_x(s):=\widetilde{q}(S|x,s)$.

Given a $Q$-function $q$ on $S$ from $S\times [0,\infty)$, for each $\Gamma\in {\cal B}(S)$, $x\in S$, $s,t\in [0,\infty)$ and $s\le t,$ one can define
\begin{eqnarray*}
&&p_{q}^{(0)}(s,x,t,\Gamma):= \delta_{x}(\Gamma) e^{-\int_s^t q_x(v)dv},\nonumber\\
&&p_q^{(n+1)}(s,x,t,\Gamma):= \int_s^t e^{-\int_s^u q_x(v)dv} \left(\int_S  p_q^{(n)}(u,z,t,\Gamma)\widetilde{q}(dz|x,u)\right)du,~\forall~ n=0,1,\dots.
\end{eqnarray*}
It is clear that one can legitimately define the sub-stochastic kernel $p_q(s,x,t,dy)$ on $S$ by
\begin{eqnarray}\label{ZY2015:01}
p_q(s,x,t,\Gamma):=\sum_{n=0}^\infty p_q^{(n)}(s,x,t,\Gamma)
\end{eqnarray}
for each $x\in S$, $s,t\in [0,\infty),$ $s\le t$, and $\Gamma\in {\cal B}(S)$. This is the Feller's construction for a transition function, i.e., $p_q$ satisfies
\begin{eqnarray*}
p_q(s,x,s,dy)=\delta_x(dy)
\end{eqnarray*}
and the Kolmogorov-Chapman equation
\begin{eqnarray*}
\int_S p_q(s,x,t,dy)p_q(t,y,u,\Gamma)=p_q(s,x,u,\Gamma),~\forall~\Gamma\in {\cal B}(S)
\end{eqnarray*}
is valid for each $0\le s\le t\le u<\infty$.

\begin{condition}\label{ZY2015Condition2}
There exist a monotone nondecreasing sequence $\{\widetilde{V}_n\}_{n=1}^\infty\subseteq {\cal B}([0,\infty)\times S)$ and a $[0,\infty)$-valued measurable function $\widetilde{w}$ on $[0,\infty)\times S$ such that the following hold.
\begin{itemize}
\item[(a)]  As $n\uparrow \infty,$ $\widetilde{V}_n\uparrow [0,\infty)\times S$.
\item[(b)] For each $n=1,2,\dots,$
$
\sup_{x\in \hat{V}_n,~t\in[0,\infty)} q_x(t)<\infty,
$
where $\hat{V}_n$ denotes the projection of $\widetilde{V}_n$ on $S$.
\item[(c)] As $n\uparrow \infty,$
$
\inf_{(t,x)\in ([0,\infty)\times S)\setminus \widetilde{V}_n} \widetilde{w}(t,x)\uparrow \infty.
$
\item[(d)] For some constant $\rho' \in (0,\infty),$ for each $x\in S$ and $v\in[0,\infty)$,
\begin{eqnarray*}
\int_0^\infty\int_{S} \widetilde{w}(t+v,y) e^{-\rho' t -\int_0^t q_x(s+v)ds} \widetilde{q}(dy|x,t+v)dt\le \widetilde{w}(v,x).
\end{eqnarray*}
\end{itemize}
\end{condition}

The next statement follows from Theorem 3.2 of \cite{Zhang:2016}.
\begin{theorem}\label{ZY2015Theorem03}
If Condition \ref{ZY2015Condition2} is satisfied, then $p_q(s,x,t,S)=1$ for each $x\in S$, $s,t\in[0,\infty)$ such that $s\le t.$
\end{theorem}

\par\noindent\textbf{Acknowledgement.} This work is partially supported by a grant from the Royal Society (IE160503).


\begin{thebibliography}{99}

\bibitem{Altman:1999} Altman, E. (1999). {\em Constrained Markov Decision Processes}. Chapman and Hall/CRC, Boca Raton.

\bibitem{Anderson:1991} Anderson, W. (1991). {\em Continuous-time Markov Chains}. Springer, New York.

\bibitem{BauerleRieder:2011}  B\"auerle, N. and Rieder, U.
(2011). {\em Markov Decision Processes with Applications to
Finance}. Springer, Berlin.

\bibitem{Blok:2015} Blok, H. and Spieksma, F. (2015). Countable state Markov decision processes with unbounded jump rates and discounted optimality equation and approximations. {\em Adv. Appl. Probab.} {\bfseries 47}, 1088-1107.

\bibitem{Chen:2015} Chen, M. (2015). Practical criterion for uniqueness of q-processes. {\em Chinese J. Appl. Probab. Stat.}, {\bfseries 31}, 213--224.

\bibitem{CostDufour:2014} Costa, O. and Dufour, F. (2015). A linear programming formulation for constrained discounted continuous control for piecewise deterministic Markov processes. {\em J. Math. Anal. Appl.} {\bfseries 424}, 892-914.


\bibitem{Dufour:2012} Dufour, F., Horiguchi, M. and Piunovskiy, A. (2012). {The expected total cost criterion for Markov decision processes under constraints: a convex analytic approach}.
{\em Adv. Appl. Probab.} {\bfseries 44}, 774--793.

\bibitem{Dufour:2013} Dufour, F. and Piunovskiy, A. (2013). The expected total cost criterion for Markov decision processes under constraints. {\em Adv. Appl. Probab.} {\bfseries 45}, 837-859.

\bibitem{Dufour:2016} Dufour, F. and Prieto-Rumean, T. (2016) Conditions for the solvability of the linear programming formulation for constrained discounted Markov decision processes. {\em Appl. Math. Optim.} {\bfseries 74}, 27-51.

\bibitem{Dynkin:1979} Dynkin, E. and Yushkevich, A. (1979). {\em Controlled Markov Processes}. Springer, New York.

\bibitem{FeinbergSonin:1996} Feinberg, E. and Sonin, I. (1996). Notes on equivalent stationary policies in Markov decision processes with total rewards. {\em Math. Meth. Oper. Res.} {\bfseries 44}, 205-221.

\bibitem{Feinberg:2004} Feinberg, E. (2004). Continuous time discounted jump Markov decision
processes: a discrete-event approach. {\em Math. Oper. Res.}
{\bfseries 29}, 492-524.


\bibitem{Feinberg:2012OHL} Feinberg, E. (2012). Reduction of discounted continuous-time MDPs with unbounded jump
and reward rates to discrete-time total-reward MDPs. In {\em
Optimization, Control, and Applications of Stochastic Systems},
Hernandez-Hernandez, D. and Minjarez-Sosa, A. (eds): 77-97,
Birkh\"{a}user, Bassel.



\bibitem{Feinberg:2012} Feinberg, E. and Rothblum, U. (2012). Splitting randomized stationary policies in total-reward Markov
decision processes. {\em Math. Oper. Res.} {\bfseries 37},
129-153.

\bibitem{FeinbergShiyaevGood:2013} Feinberg, E., Mandava, M. and Shiryaev, A. (2013). Sufficiency of Markov policies for continuous-time Markov decision processes and solutions of Kolmogorov's forward equation for jump Markov processes. In {\em Proc. 52nd IEEE CDC}, 5728-5732. Dec, 2013, Florence, Italy.

\bibitem{FeinbergShiryayev:2014} Feinberg, E., Mandava, M. and
Shiryaev, A. (2014). On solutions of Kolmogorov's equations for
nonhomogeneous jump Markov processes. {\em J. Math. Anal. Appl.},
{\bfseries 411}, 261--270.

\bibitem{Gihman:1975} Gihman, I. and Skorohod, A. (1975). {\em The Theory of Stochastic Processes II}. Springer, Berlin.

\bibitem{Guo:2007} Guo, X. (2007). Continuous-time Markov decision processes with discounted rewards: the case of Polish spaces. {\em Math. Oper. Res.}, {\bfseries32}, 73-87.



\bibitem{Guo:2009} {Guo, X. and Hern\'{a}ndez-Lerma, O.} (2009).
{\em {C}ontinuous-Time {M}arkov {D}ecision {P}rocesses: {T}heory
and {A}pplications}. Springer, Heidelberg.

\bibitem{GuoZhang:2016} Guo, X. and Zhang, Y. (2016). Constrained total undiscounted continuous-time Markov decision processes. {\em Bernoulli}, accepted.

\bibitem{Hernandez-Lerma:1996} Hern{\'a}ndez-Lerma, O. and Lasserre, J. (1996). {\em Discrete-Time {M}arkov Control
Processes}. Springer-Verlag, New York.

\bibitem{Jacod:1975} Jacod, J. (1975).
Multivariate point processes: predictable projection,
Radon-Nykodym derivatives, representation of martingales. {\em Z.
Wahrscheinlichkeitstheorie verw. Gebite.} {\bfseries 31},
235-253.

\bibitem{Jaskiewicz:2011} Ja{\'s}kiewicz, A. and Nowak, A. (2011). Stochastic games with unbounded payoffs: applications to robust control in economics. {\em Dyn. Games Appl.} {\bfseries 1}, 253-279.


\bibitem{Kitaev:1995} Kitaev, M. and Rykov, V. (1995). {\em Controlled Queueing Systems}. CRC Press, Boca Raton.

\bibitem{Kuznetsov:1984} Kuznetsov, S. (1984). Inhomogeneous Markov processes. {\em J. Soviet Math.} {\bfseries 25}, 1380--1498.



\bibitem{Piunovskiy:1997} Piunovskiy, A. (1997). {\em Optimal Control of Random Sequences in Problems with Constraints}, Kluwer, Dordrecht.

\bibitem{PiunovskiyZhang:2011Arxiv} Piunovskiy, A. and Zhang, Y.
(2011). Discounted continuous-time markov decision processes with
unbounded rates: the dynamic programming approach. Available at arXiv:1103.0134.


\bibitem{ABPZY:20102} Piunovskiy, A. and Zhang, Y. (2011). Discounted continuous-time Markov decision processes with unbounded rates: the convex analytic approach. {\em SIAM J. Control Optim.} {\bfseries 49}, 2032-2061.


\bibitem{Piunovskiy:2015} Piunovskiy, A. (2015). Randomized and relaxed strategies in continuous-time Markov decision processes. {\em SIAM J. Control Optim.} {\bfseries 53}, 3503-3533.

\bibitem{PrietoOHL:2012book}  Prieto-Rumeau, T. and
Hern{\'a}ndez-Lerma, O. (2012). {\em Selected Topics in
Continuous-Time Controlled Markov Chains and Markov Games}.
Imperial College Press, London.


\bibitem{Rykov:1966} Rykov, V. (1966). Markov decision processes with finite state and decision spaces.
{\em Theory Probab. Appl.} {\bfseries 11}, 302–-311.


\bibitem{Spieksma:2015} Spieksma, F. (2015). Countable state Markov processes:
non-explosiveness and moment function. {\em Probab. Eng. Inform. Sc.}, {\bfseries 29}, 623--637.


\bibitem{Spieksma:2012} Spieksma, F. (2016). Kolmogorov forward equation and explosiveness in countable state Markov processes. {\em Ann. Oper. Res.} {\bfseries 241}, 3--22.

\bibitem{vanderWal:1981} van der Wal, J. (1980). {\em Stochastic Dynamic Programming: Successive Approximations and Nearly Optimal Strategies for Markov Decision Processes and Markov Games}. Mathematisch Centrum, Amsterdam.

\bibitem{Zhang:2016} Zhang, Y. (2015). On the nonexplosion and explosion for nonhomogeneous Markov pure jump processes. Preprint. Available at arXiv:http://arxiv.org/abs/1511.05011.

\end{thebibliography}
\end{document}